\documentclass{amsart}
\usepackage{amssymb, amsmath}

\def\largerightarrow{   -\negthinspace\negthinspace -\negthinspace\negthinspace
                                       \negthinspace\longrightarrow }

\begin{document}
\Large
\title{}
\vbox{\hfil {\Large\bf ON RATIONAL INJECTIVITY OF KASPAROVS}\hfil}
\vbox{\hfil {\Large\bf  ASSEMBLY MAP IN DIMENSION $\,\leq 2\,$}\hfil}

\author{U. Haag}

\date{\today \\ \texttt{\hfil Contact:haag@mathematik.hu-berlin.de}}
\maketitle
\medskip\medskip\noindent
\hfil\hbox{ {\it To Joachim Cuntz, teacher and friend, on the occasion}}\hfil
\par
\hfil\hbox{ {\it of his sixtieth birthday, and to Gennadi Kasparov, whose}}\hfil
\par
\hfil\hbox{ {\it work is a great source of inspiration. }}\hfil
\bigskip\bigskip\noindent
\begin{abstract}
In this article the author presents a new proof of injectivity of the composition of the inverse of the rational Chern  Character in homology applied to the classifying space $\, BG\, $ of a (countable) discrete group $\, G\, $ and restricted to dimensions less or equal than two with the rationalized assembly map of Kasparov, (cf. \cite{K2})
$$ H_{\leq 2} ( BG ; {\mathbb Q} ) \buildrel ch^{-1}_{\mathbb Q}\over\largerightarrow 
RK_* ( BG ) \otimes {\mathbb Q} \buildrel A\over\largerightarrow 
K_* ( C^* ( G ) ) \otimes {\mathbb Q}\, . $$
\par\medskip\noindent
Here  the lefthandside is the assemblage of its ordinary homology groups with rational coefficients of dimensions $\,\leq 2\, $, whereas the middle term is given by  the direct limit of the 
analytic $K$-homology groups of Kasparov taken with respect to the 
compact subsets  of $\, BG\, $ (ordered by inclusion) and tensored with the rationals, and 
the righthandside denotes the (operator) $K$-theory of the full group $C^*$-algebra of 
$\, G\, $, again tensored by the rational numbers.
\end{abstract}
\par\bigskip\bigskip\bigskip\noindent
{\bf Foreword.}\quad This article was intended to appear in the special edition of the Mathematical 
Journal of the University of Muenster, Westfalen, in honour of Joachim Cuntz, but at that time the main argument still contained some severe defects which could not be overcome easily, so the paper had to be withdrawn. Some time has elapsed since then and the author hopes that the arguments presented here are correct and understandable to the interested reader. Despite the efforts for clarity that have been made the author must admit, and regrets, that some of the references given are not easily available. This concerns the paper "KK-theory of stable projective limits" of which the first part-
dealing with the covariant side of the KK-functor- has appeared in the preprint series of the IML at 
Luminy, Marseille. The second part dealing with K-homology only exists as a manuscript of the author, but anyway is not relevant for the purposes of this paper which only uses results on the covariant side. 
It is planned to submit the complete paper for publication in the near future. Similarly, the other reference 
by the author, the treatise on "Regular algebraic K-theory for groups", is unpublished up 
to this date. Only some very minor results of this paper are used -- all of them dealing with low dimensions where the algebraic K-theory groups coincide with ordinary group homology -- so that a standard textbook on Homological Algebra (cf. \cite{Hi}) giving the basic definitions of group homology including exact sequences in  low dimensions associated with group extensions might be helpful. In addition, the interested reader is encouraged to contact the author for copies of the articles in question.  
\par\bigskip\noindent
{\bf 0. Introduction}
\par\smallskip\noindent
There is by now an overwhelming variety of special classes of groups for which the (integral) injectivity of the assembly map as above has been proved in full generality without any restrictions to certain dimensions and in many cases even stronger results are available 
(Baum-Connes type isomorphism theorems). On the other hand up to this date and knowledge of the author not a single counterexample has been constructed where the full integral injectivity of the assembly map (replacing $\, BG\, $ with the classifying space for proper actions in case of a group with torsion) fails. The proofs however all seem to depend on certain geometric properties of the group 
$\, G\, $ and it is questionable whether this method will eventually succeed in covering all 
(say finitely generated discrete) groups.  The Novikov Conjecture on homotopy invariance of higher signatures only depends on rational injectivity of the assembly map which is known for general groups in low dimensions. The first step was taken by Novikov who observed the homotopy invariance of higher signatures associated with onedimensional cohomology classes  (the zerodimensional case being trivial), followed by Kasparov who extended the result to arbitrary products of onedimensional classes. The twodimensional case appears in Mathai \cite{Ma} and a paper of
Hanke and Schick \cite{H-S} where also the injectivity of the rational assembly map restricted to elements which are dual to the subring of cohomology generated by the one- and twodimensional classes is established. 
The present paper gives a new proof of a somewhat more modest result, namely rational injectivity
for twodimensional homology classes (and duals of products of onedimensional cohomology classes).   
It should be pointed out that this is not sufficient to prove homotopy invariance of the corresponding higher signatures which would however follow from a dual result (rational surjectivity of the dual assembly map) which is fairly easy to obtain in dimension one as sketched below, thereby 
confirming Kasparovs result on such classes.
\par\bigskip\noindent
{\bf 1. Outline of the argument}
\par\smallskip\noindent
We begin by recalling some results of Appendix C of the paper \cite{Ha1} by the author. In this paper there is defined for any normal subgroup $\, N\, $ of a discrete group $\, G\, $ a sequence of abelian groups 
$\{ K^J_n ( N , G ) {\} }_{n \geq 1} $ called the (relative) $K^J$-groups of the pair $\, ( N , G )\, $ and behaving functorially with respect to (compatible) group homomorphisms. In particular one has for any group 
$\, G\, $ the absolute $K^J$-groups $\, K^J_n ( G ) := K^J_n ( G , G )\, $. To any group extension 
$$ 1 \largerightarrow N \largerightarrow G \largerightarrow H \largerightarrow 1 $$
there is associated a long sequence of $K^J$-groups of the form
$$ \cdots \longrightarrow K^J_{n+1} ( H ) \buildrel {\partial }_n\over\longrightarrow K^J_n ( N , G ) \longrightarrow K^J_n ( G ) \longrightarrow K^J_n ( H ) \buildrel {\partial }_{n-1}\over\longrightarrow \cdots $$
which is everywhere exact in certain cases (for example if $\, G\, $ and hence $\, H\, $ are perfect groups and $\, N\, $ is {\it full}, i.e. equal to the commutator subgroup $\, [ N , G ]\, $ generated by commutators 
$\, [ x , y ] = x\, y\, x^{-1}\, y^{-1}\, $ with $\, x\in N\, ,\, y\in G\, $). For a general (normal sub)group the sequence is also everywhere exact except possibly at 
$\, K^J_3 ( N , G )\, $ and $\, K^J_3 ( G )\, $. There is a natural transformation to group homology which in case of dimensions one and two reduces to a canonical isomorphism. In particular one gets that 
$\, K^J_2 ( G ) = H_2 ( G ) = H_2 ( BG )\, $ is given by the Schur multiplier of the group $\, G\, $. The general theory and possible failure of exactness for general groups explained above is but of no importance to the present argument since the sequence will be exact in any case in dimensions one and two with which we are dealing here. What matters is that one constructs a canonical (mod 2-) degreepreserving map into $\, RK_* ( BG )\, $ for each group $\, K^J_n ( G )\, $. One easily finds that the image of $\, K^J_2 ( G )\, $ is rationally just the preimage of $\, H_2 ( BG )\, $ for the Chern Character. By construction this map commutes with the boundary maps in the following sense:
given a group extension as above one obtains a sixterm exact sequence in $RK$-homology which relates the corresponding groups for $\, BG\, $, $\, BH\, $ and the mapping cone $\, C_{\varphi }\, $ of the induced map $\, \varphi : BG \rightarrow BH\, $ (restricting to compact subsets and passing to the $C^*$-algebras of continuous functions etc.). Then one has a commutative diagram
$$ \vbox{\halign{ #&#&#\cr
 \hfil $K^J_{n+1} ( H )$\hfil &\hfil $\buildrel {\partial }_n\over\largerightarrow $\hfil 
 &\hfil $K^J_n ( N , G )$\hfil \cr
 \hfil $\Bigm\downarrow $\hfil &&\hfil $\Bigm\downarrow $\hfil \cr
 \hfil $RK_{*+1} ( BH )$ \hfil &\hfil $\buildrelÊ{\partialÊ}_*\over\largerightarrow $\hfil &
 \hfil $RK_* ( C_{\varphi } )$ \hfil \cr  }} . $$
 \par\medskip\noindent
 Composing this canonical map with the (integral) assembly map results in a natural transformation from 
 $\, K^J_{2n+*} ( G )\, $ to $\, K_* ( C^* ( G ))\, $ (note that both 
 $\, BG\, $ and $\, C^* ( G )\, $ behave completely functorial with respect to group homomorphisms). Again this will result in a natural map from 
 $\, K^J_n ( N , G )\, $ to the $K$-group of proper dimension of the mapping cone of the induced homomorphism $\, C^* ( G ) \twoheadrightarrow C^* ( H )\, $ which for the sake of $K$-theory is equivalent to its kernel. Consider a free resolution 
 $$ 1 \longrightarrow R \longrightarrow F \longrightarrow G \longrightarrow 1 $$
 of $\, G\, $ and put $\, {\overline G} := F / [ R , F ]\, $ where $\, [ R , F ]\, $ is the commutator subgroup generated by all single commutators $\, r\, x\, r^{-1}\, x^{-1}\, $ with $\, r \in R\, ,\, x \in F\, $.  Clearly this is a central extension of $\, G\, $ whose kernel $\, K\, $ contains $\, K^J_2 ( G )\, $ as a direct summand complemented by a free abelian group. (Recall that the Schur multiplier is given by the formula 
$$ K^J_2 ( G ) := R \cap [ F , F ] / [ R , F ] \> .\> )$$
Considering the long sequence of $K^J$-theory applied to the extension
$$ 1 \longrightarrow K \longrightarrow {\overline G} \longrightarrow G \longrightarrow 1 $$
one gets that $\, K^J_1 ( K , {\overline G} ) := K / [ K , {\overline G} ] = K\, $ because $\, K\, $ is central and the boundary map $\, K^J_2 ( G ) \buildrel\partial\over\longrightarrow K^J_1 ( K , {\overline G} )\, $ is an injection. Let $\, C_0^* ( G )\, $ denote the augmentation ideal in $\, C^* ( G )\, $, i.e. the kernel of the natural map $\, C^* ( G ) \twoheadrightarrow {\mathbb C}\, $ induced from the trivial representation of 
$\, G\, $ so that the kernel of $\, C_0^* ( {\overline G} ) \twoheadrightarrow C_0^* ( G )\, $ can be naturally identified with the kernel of the quotient map on the full $C^*$-algebras. It is shown in \cite{Ha1} that the image of the maps
 $\, K^J_n ( G ) \longrightarrow K_* ( C^* ( G ))\, $ actually lie in the splitting subgroup 
 $\, K_* ( C_0^* ( G ))\, $ for $\, n \geq 1\, $. Let $\, {\mathcal J}\, $ denote the kernel of the extension of group $C^*$-algebras as above and consider the commutative diagram 
$$ \vbox{\halign{ #&#&#\cr
\hfil $K^J_1 ( K )$\hfil &\hfil $\largerightarrow $\hfil &\hfil $K_1 ( C_0^* ( K ) )$\hfil \cr
\hfil $\Bigm\Vert $\hfil &&\hfil $\Bigm\downarrow $\hfil \cr
\hfil $K^J_1 ( K , {\overline G} )$\hfil &\hfil $\largerightarrow $\hfil &\hfil $K_1 ( \mathcal J )$\hfil \cr }} . $$
\par\medskip\noindent
The upper horizontal map is known to be rationally injective since $\, K\, $ is an abelian group. Thus if we can show that the composition of this map with the right vertical map is still rationally injective then by a commutative diagram considered above intertwining the boundary maps one gets that also the map 
$\, K^J_2 ( G ) \longrightarrow K_0 ( C_0^* ( G ) )\, $ must be rationally injective.
\par\smallskip\noindent
Before entering into the proof let us consider the onedimensional case which is much simpler. In fact by functoriality one gets a commutative diagram 
$$ \vbox{\halign{ #&#&#\cr
\hfil $K^J_1 ( G )$\hfil &\hfil $\largerightarrow $\hfil &\hfil $K_1 ( C^* ( G) )$\hfil \cr
\hfil $\Bigm\Vert $\hfil &&\hfil $\Bigm\downarrow $\hfil \cr
\hfil $K^J_1 ( G / [ G , G ] )$\hfil &\hfil $\largerightarrow $\hfil &\hfil $K_1 ( C^* ( G / [ G , G ] ) )$\hfil \cr }}  $$
\par\medskip\noindent
since $\, K^J_1 ( G ) = G / [ G , G ] = K^J_1 ( G / [ G , G ] )\, $ and the lower horizontal map is rationally injective as $\, G / [ G , G ]\, $ is abelian. Considering the dual assembly map which in general takes the form
$$ K^* ( C^* ( G ) ) \buildrel \widehat A\over\largerightarrow \lim_{\leftarrow X}  K^* ( X ) $$
with $\,  X \subseteq BG \, $ ranging over the compact subsets one rationally has 
$\, \displaystyle\lim_{\leftarrow X}\,  ( K^* ( X ) \otimes {\mathbb Q}) \simeq \lim_{\leftarrow X}\, 
( \oplus_n\, H^n ( X ; {\mathbb Q} ) )\, $ and projecting the righthandside onto the onedimensional 
cohomology subgroup gives the dual space of $\, H_1 ( BG ; {\mathbb Q} )\, $ (at least when 
$\, G\, $ is finitely generated) by the Universal Coefficient Theorem. Since for finitely generated abelian groups the rational dual assembly map is known to be an isomorphism and 
in general $\, H^1 ( G ; {\mathbb Q} ) = H^1 ( G / [ G , G ] ; {\mathbb Q} )\, $ one sees that this subgroup is in the image of the dual assembly map. Therefore the pairing between homology and arbitrary products of onedimensional cohomology classes can be pulled back to a pairing of 
$\, K_* ( C^* ( G ) ) \otimes {\mathbb Q}\, $ and certain classes in 
$\, K^* ( C^* ( G ) ) \otimes {\mathbb Q}\, $ which gives Kasparovs result mentioned in the Introduction. 
\par\bigskip\noindent
{\bf 2. The proof}
\par\smallskip\noindent 
In view of the discussion above we need to prove the following result
\par\medskip\noindent
{\bf Proposition.}\quad The natural composition
$$ K  \largerightarrow K_1 ( C_0^* ( K ) ) \largerightarrow K_1 ( \mathcal J ) $$
is rationally injective.
\par\bigskip\noindent 
{\it Remark.}\quad It can be shown that the map 
$\, K_1 ( C_0^* ( K ) ) \longrightarrow K_1 ( \mathcal J )\, $ is not in general rationally injective (for example on the image of $\, K^J_3 ( K )\, $).
\par\medskip\noindent
Let $\, H\, $ be any central extension of $\, G\, $ by a (countable) abelian group $\, K\, $ so that 
$\, C^* ( H )\, $ has the structure of a semicontinuous field of $C^*$-algebras $\, A\, $ over the compact dual $\,\widehat K\,$,  
where semicontinuous field simply means that $\,A\, $ is a $C ( \widehat K )$-algebra, which is the maximal $C^*$-completion of the algebra of finite sums  
$\,\{\, \sum_{g \in G}\, a_g\, s_g\,\}\, $ with coefficients in the central subalgebra 
$\, C ( \widehat K )\, $ of continuous functions on $\,\widehat K\, $ and $\, s_g\, $ is the unitary generator of $\, C^* ( H )\, $ corresponding to the group element $\, g \in G\, $ by  some chosen section 
$\, s : G \nearrow H\, $ to the natural projection. Denote the fibres of the semicontinuous field   
by $\,\{ A_{\omega } {\} }_{\omega \in \widehat K}\, $. Then the evaluated elements 
$\, \{ s_g^{\omega } \}\, $ are unitary generators of $\, A_{\omega }\, $, and for each pair of points 
$\, \omega\, ,\, {\omega }' \in \widehat K\, $ there is a natural densely defined map 
$\, A_{\omega } \rightarrow A_{{\omega }'}\, $ given by linear extension of the assignment 
$\, s_g^{\omega } \mapsto s_g^{{\omega }'}\, $ with domain the finite linear combinations 
$\, \sum_g\, a_g\, s_g^{\omega }\, $ which will be called the tautological shift.
Also on $\, A\, $ there is a $C ( \widehat K )$-valued trace $\, Tr\, $ of norm $1$ given by evaluating a finite sum $\, \sum_g\, a_g\, s_g\, $ at $\, a_1\, $ (and extending by continuity). It is well known that the norm-function $\, \omega \mapsto \Vert x {\Vert }_{\omega }\, $ is upper semicontinuous 
for any element $\, x\, $ of the semicontinuous field(cf. \cite{Bl}, \cite{Rf} ).
\par\bigskip\noindent
{\it Proof of Proposition.}\quad Let us first deal with the case where $\, K\, $ is finitely generated which is a bit more simple. Assume that the map is not rationally injective. Then there exists a nontorsion element of $\, K\, $ which is in the kernel of the composition as above. Without loss of generality we may divide 
$\, K\, $ by its torsion subgroup, and since $\, K\, $ is finitely generated we may assume that our chosen element is an element of a basis of the remaining free abelian quotient. Thus completing to a full basis we may further divide by a complementary subgroup obtaining a central extension of $\, G\, $ by 
$\, \mathbb Z\, $. We continue to denote this extension by $\, \overline G\, $ and its kernel by $\, K\, $. If the Proposition is valid for the new extension, it is also valid for the original one since the chosen 
element was completely deliberate. Now $\, C^* ( K )\, $ can be identified with continuous functions on the circle $\, S^1 = \{ z \in \mathbb C\, \vert \, \vert z\vert = 1 \}\, $ and $\, C_0^* ( K )\, $ corresponds to those functions vanishing at the point $1\, $.
The kernel $\, \mathcal J\, $ is a certain $C^*$-completion of the algebra of finite sums 
$\,\{\, \sum_{g \in G}\, a_g\, s_g\,\}\, $ with coefficients $\, a_g\, $ taken from the central subalgebra 
$\, C_0^* ( K )\, $ and $\, s_g\, $ is the unitary generator of $\, C^* ( \overline G )\, $ corresponding to the group element $\, g \in G\, $ by  some chosen section $\, s : G \nearrow \overline G\, $ to the natural projection. Multiplication depends on the particular section, i.e. for each pair 
$\, g , h \in G\, $ there is an element $\, n ( g , h ) \in \mathbb Z\, $ such that 
$\, s_g \cdot s_h = u^{n ( g , h )}\, s_{gh}\, $ with $\, u\, $ the canonical unitary generator of 
$\, C^* ( K )\, $ acting as a multiplier on $\, C_0^* ( K )\, $. One may continuously deform this action with respect to some parameter $\, \theta \in [ 0 , 1 ]\, $ by letting $\, u_{\thetaÊ}\, $ denote the multiplier of 
$\, C_0 ( S^1\backslash \{ 1 \} ) \simeq C_0^* ( K )\, $ such that $\, u_{\theta } \equiv 1\, $ for all points 
$\, z = e^{2\pi i \varphi}\, $ with $\, \varphi \leq \theta\, $ and rotating the identical function $\, u\, $ by an angle of $\, 2\pi i \theta\, $ for all other values $\,\geq\theta\, $. One gets a continuous deformation with 
$\, u_0 = u\, $ and $\, u_1 \equiv 1\, $. Denote the corresponding semicontinuous field of $C^*$-algebras by 
$\,\{ {\mathcal J}_{\theta }\}\, $. Then $\, {\mathcal J}_1\, $ is just the kernel of 
$\, C^* ( K \times G ) \twoheadrightarrow C^* ( G )\, $ and by the splitting projection 
$\, K \times G \twoheadrightarrow K\, $ the map $\, C_0 ( K ) \longrightarrow {\mathcal J}_1\, $ is seen to split, whence the corresponding map in $K$-theory is (split) injective. Now the $K$-theory of a 
semicontinuous field of $C^*$-algebras posesses some semicontinuity properties. 
If an element of 
$\, K_* ( C_0^* ( K ) )\, $ is in the kernel of the map to $\, K_* ( {\mathcal J}_{{\theta }_0} )\, $ for some fixed value $\, {\theta }_0\, $, then there exists a neighbourhood of $\, {\theta }_0\, $ such that the element is in the kernel of the map to $\, K_* ( {\mathcal J}_{\theta } )\, $ for every $\, \theta \, $ in this neighbourhood. This is seen considering the associated semicontinuous field of mapping cones 
$\, C_{\theta }\, $ for the inclusions $\, C_0^* ( K ) \hookrightarrow {\mathcal J}_{\theta }\, $. One gets surjective maps $\, C_{\theta } \twoheadrightarrow C_0^* ( K )\, $ and the question whether an element is in the kernel of the map on $K$-groups as above for some given parameter boils down up to stabilization etc. to the fact whether a given unitary (resp. projection or difference of projections) in the quotient lifts to a unitary (projection) in the corresponding mapping cone extension. But this lifting property will extend to some neighbourhood of $\, {\theta }_0\, $ from functional calculus. 
From Bott Periodicity it suffices to do the case of $\, K_1\, $. Suppose a unitary $\, u\, $ lifts to a unitary (invertible) element at some given point $\, {\theta }_0\, $. 
Consider an arbitrary lift $\, V\, $ which is invertible at $\, {\theta }_0\, $. Then $\, V_{\theta }\, $ is invertible iff  $\,\vert V{\vert }_{\theta }\, $ and $\, \vert V^* {\vert }_{\theta }\, $ are both invertible. Suppose that in each neighbourhood of 
$\, {\theta }_0\, $ there exist parameters $\,\theta \, $ where $\, \vert V {\vert }_{\theta }\, $ (or 
$\, \vert V^* {\vert }_{\theta }\, $) is not invertible.
Applying functional calculus to $\, \vert V\vert \, $ with respect to a continuous function $\, f\, $ vanishing outside the open unit disc $\, \{ z \in\mathbb C\,\vert\, \vert z\vert < 1 \}\, $ and such that $\, f \equiv 1\, $ in some neighbourhood of $\, 0\, $ one gets a contradiction to upper semicontinuity of the norm-function of 
$\, f ( \vert V\vert )\, $ at $\, {\theta }_0\, $ since this element is zero evaluated at $\,{\theta }_0\, $. 
But then there exists a unitary $\, U_{\theta }\, $ (using polar decomposition of $\, V\, $ in a given neighbourhood where it is invertible) for each point in the neighbourhood such that its image maps to the image of $\, u\, $ in the
$K_1$-group of the quotient evaluated at the given point. 
Thus one finds that the set of parameters $\, \theta\, $ where the map on $K$-groups is injective is closed, since its complement is open. We are interested in injectivity for the parameter $\, \theta = 0\, $, so it suffices to show that the map is injective for all positive parameters. Obviously for every positive parameter $\, \theta\, $ there is an interval on the circle where the action of the multiplier $\, u_{\theta }\, $ reduces to the identity, moreover we can easily contract the original embedding 
$\, C_0^* ( K ) \hookrightarrow {\mathcal J}_{\theta }\, $ by homotopy to an embedding which is trivial outside the interval in question. However we cannot just cut off the complementary interval so as to reduce to the case of the parameter $\, \theta = 1\, $. Put 
$\, B = C_0^* ( K )\, $. We will now show the converse, namely that if the map $\, K_1 ( B ) \rightarrowtail K_1 ( {\mathcal J}_{{\theta }_0} )\, $ is injective for some parameter $\, {\theta }_0\in [ 0 , 1Ê]\, $ then there exists a neighbourhood $\, {\mathcal U}_{{\theta }_0}\, $ of $\,{\theta }_0\, $ where the map $\, K_1 ( B ) \rightarrowtail K_1 ( {\mathcal J}_{\theta } )\, $ is injective for all $\, \theta \in {\mathcal U}_{{\theta }_0}\, $. To this end assume given 
$\, {\theta }_0\, $ such that the first map as above is injective. Clearly this implies that the map 
$\, K_1 ( B ) \rightarrowtail K_1 ( {\mathcal J}_{\theta } ) \, $ is also injective for each 
$\, \theta \geq {\theta }_0\, $ because the first map will factor over the second one by the embedding of 
$\, {\mathcal J}_{\theta }\, $ as an ideal of $\, {\mathcal J}_{{\theta }_0}\, $ using some convenient reparametrization. Now assume that injectivity fails for every $\, \theta < {\theta }_0\, $. 
Choose a monotonously increasing sequence $\,\{ {\theta }_k {\} }_{k\geq1}\, $ of points converging up to 
$\, {\theta }_0\, $.
The assumption implies on embedding $\, {\mathcal J}_{{\theta }_0}\, $ as an ideal in 
$\, {\mathcal J}_{\theta }\, $ and using the six-term exact sequence of $K$-theory, that the generator of the kernel of the composition 
$$K_1 ( B ) \rightarrowtail K_1 ( {\mathcal J}_{{\theta }_0} ) \rightarrow K_1 ( {\mathcal J}_{\theta } )$$
has a lift to the group $\, K_0 ( {\mathcal J}_{\theta } / {\mathcal J}_{{\theta }_0} )\, $ for any 
$\, \theta < {\theta }_0\, $. This quotient algebra is again a semicontinuous field over the halfopen interval 
$\, ( \theta\, ,\, {\theta }_0 ]\, $. On tensoring with $\, \mathcal K\, $ we can assume that all algebras are stable so that the set of all 
$\, \mathcal A_n = \mathcal K \otimes ({\mathcal J}_{{\theta }_1} / {\mathcal J}_{{\theta }_n})\, $ defines a stable projective system of $C^*$-algebras in the sense of \cite{Ha2} and there is a natural map from 
$\, \mathcal A_0 = \mathcal K \otimes ({\mathcal J}_{{\theta }_1} / {\mathcal J}_{{\theta }_0})\, $ into its stable projective limit
$\, \mathcal A_\infty = \displaystyle{ \lim_{\leftarrow }}\, \mathcal A_n\, $. Put 
$\, {\mathcal S}_n = {\mathcal J}_{{\theta }_n}\, $ and 
$\, {\mathcal C}_n = {\mathcal S}_n / {\mathcal S}_0\, $ which can be viewed as semicontinuous fields over the open interval $\, ( {\theta }_n , 1 + {\theta }_n )\, $ resp. the halfopen interval 
$\, ( {\theta }_n , {\theta }_0 ]\, $ (the parametrization is adjusted to the multiplier action of 
$\, u\, $ on the fibres of $\,\mathcal J\, $, i.e. $\, u\, $ corresponds to multiplication by the complex number 
$\, \exp (2 \pi i ( 1 - \theta ))\, $ for a parameter $\, 0 \leq \theta < 1\, $ and to the identity for 
$\, \theta \geq 1\, $). If $\, A_{\theta }\, $ denotes the fibre of 
$\,\mathcal J\, $ at the point $\, 1 - \theta\, $ there is,
for each $n$, a surjective evaluation map from $\, {\mathcal C}_n\, $ onto $\, A_{{\theta }_0}\, $ commuting with the natural inclusion maps. Let $\, \mathcal O_n\, $ be the kernel of this map. To clarify the situation consider the following three inverse systems of short exact sequences of $C^*$-algebras
$$ 1\> \longrightarrow\> \mathcal K \otimes  {\mathcal S}_n\> \longrightarrow\> 
\mathcal K \otimes {\mathcal S}_1\> \longrightarrow\> {\mathcal A}_n\> \longrightarrow\> 1 $$
$$ 1\> \largerightarrow\> {\mathcal S}_0\> \largerightarrow\> {\mathcal S}_n\> \largerightarrow\> 
{\mathcal C}_n\> \largerightarrow\> 1 $$
$$ 1\>\largerightarrow\> \mathcal O_n\>\largerightarrow\> \mathcal C_n\>\largerightarrow\> 
A_{{\theta }_0}\>\largerightarrow\> 1 \> . $$
Recall the following facts: computing the $KK$-theory (or $E$-theory) of a stable projective limit in the second 
(covariant) variable no additional assumptions on the inverse system of $C^*$-algebras 
$\, \{ A_n , {\rho }_n \}\, $ like separability or nuclearity are needed, only stability and that the connecting homomorphisms are surjective. On the other hand for any injective inverse system as given above by the $\, \{ {\mathcal S}_n \}\, $ or $\, \{ {\mathcal C}_nÊ\}\, $ such that the smaller algebras sit as ideals of the larger ones one constructs an equivalent stable projective system considering the associated system of mapping cones $\, \{ C_{{\varphi }_n} \}\, $ with 
$\, {\varphi }_n : {\mathcal K} \otimes {\mathcal S}_1 \twoheadrightarrow  
{\mathcal K} \otimes {\mathcal S}_1 / {\mathcal K}Ê\otimes {\mathcal S}_n\, $ etc.. Here equivalent means that there is a natural compatible $E$-equivalence of the factors. When alluding to the $K_*$-theory of the stable projective limit of an injective system of ideals $\,\{ \mathcal I_n \}\, $ we tacitly assume that the projective limit of the corresponding system of mapping cones is understood which is then well defined up to $E$-equivalence and denoted $\, \mathcal I_\infty\, $. From Theorem 2 of \cite{Ha2} there is for any stable projective system 
$\,\{ A_n , {\rho }_n \}\, $ a Milnor $\lim^1$-sequence of the form
$$ 1\> \longrightarrow\> {\lim}^1\, K_{*+1} ( A_n )\> \longrightarrow\> 
K_* ( \lim_{\leftarrow }\, A_n )\> \longrightarrow\> \lim_{\leftarrow }\, K_* ( A_n )\> 
\longrightarrow\> 1 \> . $$
The $\displaystyle {\lim_{\leftarrow }}$- resp $\lim^1$-groups appear as kernel resp. cokernel of the map 
$$ \prod_{n=1}^\infty\, K_* ( A_n )\> \buildrel \prod ( id - {\rho }_n )_*\over\largerightarrow\> 
\prod_{n=1}^\infty\, K_* ( A_n ) \> . $$
Recall that given a sequence of (stable) $C^*$-algebras $\,\{ A_n \}\, $ the $C^*$-algebra 
$\, {\prodÊ}_{n=1}^{\infty }\, A_n\, $ consists of norm-bounded sequences $\, x = \{ x_n \}\, $ with 
$\, x_n \in A_n\, $ and product type involution, multiplication and 
$\, \Vert x \Vert = \sup_{n\in\mathbb N} \{ \Vert x_n \Vert \}\, $. For this particular case the Milnor 
$\lim^1$-sequence yields a natural identification $\, K_* ( \Pi\, A_n )\, =\, \Pi\, K_* ( A_n )\, $.  
The $C^*$-algebra $\, \Pi\, A_n\, $ contains the direct sum $\, \oplus\, A_n\, $ as a proper ideal and the quotient will be denoted by $\, \Pi / \oplus\, A_n\, $, i.e. one has an exact sequence of the form 
$$ 1 \largerightarrow \bigoplus \, A_n \largerightarrow \prod \, A_n\largerightarrow  \prod / \bigoplus\, A_n \largerightarrow 1 \> . $$
Our first objective is to show that the natural map 
$\, K_* ( \Pi / \oplus\, {\mathcal C}_n ) \rightarrow K_* ( \Pi / \oplus\, A_{{\theta }_0} )\, $  
obtained from the evaluation map at $\, {\theta }_0\, $ is trivial. 
Since the map 
$$ K_* ( \oplus\, \mathcal C_n ) \buildrel \oplus ( id - {\rho }_n )_*\over\largerightarrow 
K_* ( \oplus\, \mathcal C_n ) $$
is an isomorphism, the kernel and cokernel of 
$$ K_* ( \Pi\, \mathcal C_n ) \buildrel \Pi ( id - {\rho }_n )_*\over\largerightarrow K_* ( \Pi\, \mathcal C_n ) $$
are the same as kernel and cokernel of 
$$ K_* ( \Pi / \oplus\, \mathcal C_n ) \buildrel \Pi / \oplus ( id - {\rho }_n )_*\over\largerightarrow 
K_* ( \Pi / \oplus\, \mathcal C_n )\> . $$
Then our assertion implies that also the evaluation map $\, K_* ( {\mathcal C}_{\infty } ) \rightarrow 
K_* ( A_{{\theta }_0} )\, $ is trivial since it factors over $\, K_* ( \Pi / \oplus\, {\mathcal C}_n )\, $ by the embedding of $\, K_* ( A_{{\theta }_0} )\, $ into $\, \Pi / \oplus\, K_* ( A_{{\theta }_0} )\, $ as the image of the constant sequences. One uses the following two facts. First, there is to each pair of indices 
$\, ( \theta\, ,\, \sigma )\, $ a *-homomorphism 
$\, {\delta }_{\sigma } : A_{\theta } \rightarrow A_{\sigma } \otimes A_{1 - \sigma + \theta }\, $ (say maximal tensor product on the right) defined by sending the unitary generators $\, s^{\theta }_g\, $ to the diagonal elementary tensors  
$\, s_g^{\sigma }\otimes s_g^{1 + \theta - \sigma }\, $, extending linearly and continuously. Second, putting $\, {\theta }' = 1-\theta\, $ there is for each index $\,\theta\, $ an antilinear isomorphism 
$$ j_{\theta } : A_{\theta } \longrightarrow A_{{\theta }'} $$
from antilinear extension of the assignment $\, s_g^{\theta } \mapsto s_g^{{\theta }'}\, $ which gives a 
$K$-equivalence of the underlying real Banach algebras and a 
$KK$-equivalence ($*$-isomorphism) of the complex $C^*$-algebra $\, A_{\theta } {\otimes }_{\mathbb R} \mathbb C\, $ with the complex $C^*$-algebra $\, A_{{\theta }'} {\otimes }_{\mathbb R} \mathbb C\, $. 
One notes that in case of a group $C^*$-algebra one has an isomorphism 
$\, C^* ( G ) {\otimes }_{\mathbb R} \mathbb C \simeq C^* ( G ) \oplus C^* ( G )\, $ which is compatible 
with the projections induced by $\, \overline G \twoheadrightarrow G\, $ so it passes to the kernel 
$\, \mathcal J {\otimes }_{\mathbb R} \mathbb C \simeq \mathcal J \oplus \mathcal J\, $. If 
$\, \tau\, $ denotes the flip of the two copies of the original group $C^*$-algebra in 
$\, C^* ( \overline G ) {\otimes }_{\mathbb R} \mathbb C\, $ then 
$\  \bigl( C^* ( \overline G ) {\otimes }_{\mathbb R} \mathbb C \bigr) {\rtimes }_{\tau } {\mathbb Z}_2 
\simeq M_2 ( C^* ( \overline G ) )\, $ and similarly with $\, \mathcal J\, $.  One has to keep in mind however that the flip automorphism is not in general compatible with evaluation at the fibres, only if one considers intervals which are symmetric around the origin.
To obtain the analogue of 
$\, {\delta }_{\sigma }\, $ in this context one uses the identity 
$\, C^* ( \overline G ) {\otimes }_{\mathbb R} \mathbb C \simeq 
C^* ( \overline G ) \oplus C^* ( \overline G )\, $. The comultiplication $\, \delta\, $ on the group 
$C^*$-algebra gives a map 
$$ C^* ( \overline G ) {\otimes }_{\mathbb R} \mathbb C \buildrel \delta \oplus \delta\over\largerightarrow 
\bigl(\, C^* ( \overline G ) \otimes C^* ( \overline G ) \,\bigr) \oplus \bigl(\, C^* ( \overline G ) \otimes 
C^* ( \overline G )\, \bigr) $$
$$ \qquad\qquad \subseteq \bigl(\, C^* ( \overline G ) {\otimes }_{\mathbb R} \mathbb C\,\bigr) \otimes \bigl(\, C^* ( \overline G ) {\otimes }_{\mathbb R} \mathbb C\,\bigr) \> .$$
Composing with the evaluation map at a point $\,\sigma\, $ in the left factor gives 
$$ C^* ( \overline G ) {\otimes }_{\mathbb R} \mathbb C \buildrel {\delta }_{\sigma }\over\largerightarrow 
\bigl(\, A_{\sigma } {\otimes }_{\mathbb R} \mathbb C \,\bigr) \otimes \bigl(\, C^* ( \overline G ) {\otimes }_{\mathbb R} \mathbb C \,\bigr) \> . $$ 
If $\, {\mathcal J}_{\theta }\, $ is an ideal of $\, C^* ( \overline G )\, $ consisting of functions vanishing at a given point $\, \theta\, $ then restriction of $\, {\delta }_{\sigma }\, $ to 
$\, {\mathcal J}_{\theta } {\otimes }_{\mathbb R} \mathbb C\, $ results in a homomorphism 
$$ {\mathcal J}_{\theta } {\otimes }_{\mathbb R} \mathbb C \buildrel {\delta }_{\sigma }\over\largerightarrow \bigl(\, A_{\sigma } {\otimes }_{\mathbb R} \mathbb C \, \bigr) \otimes 
\bigl(\, {\mathcal J}_{\theta - \sigmaÊ} {\otimes }_{\mathbb R} \mathbb C\,\bigr) $$
where the second factor consists of functions vanishing at the point $\, \theta - \sigma \, $ (or 
$\, 1 + \theta - \sigma\, $ for that matter), so $\, {\delta }_{\sigma }\, $ passes to the quotient algebras
$$ A_{\theta } {\otimes }_{\mathbb R} \mathbb C \buildrel {\delta }_{\sigma }\over\largerightarrow 
\bigl(\, A_{\sigma } {\otimes }_{\mathbb R} \mathbb C\,\bigr) \otimes \bigl(\, A_{1+ \theta - \sigma} 
{\otimes }_{\mathbb R} \mathbb C\, \bigr)\> . $$
Without loss of generality we may replace every $C^*$-algebra $\, \mathcal A\, $ considered above by 
$\, {\mathcal A} {\otimes }_{\mathbb R} \mathbb C\, $ which we will do without changing the notation. If we wish to refer to the original algebra $\, \mathcal A\, $ we use the notation $\, {\mathcal A}^0\, $. 
Now there is
a complex isomorphism $\, j_{\theta } : A_{\theta } \rightarrow A_{{\theta }'}\, $. 
All the algebras we consider may be indexed by subsets of the unit circle, the fibres $\, A_{\theta }\, $ corresponding to points. If 
$\, \mathcal A\, $ is any such algebra let $\, {\mathcal A}^{-\theta }\, $ denote the corresponding algebra with index set shifted by $\, - \theta\, $. 
Let $\, \{ A_{\lambda } {\} }_{\lambda \in \Lambda }\, $ be a cofinal net of separable $C^*$-algebras of 
$\, \Pi / \oplus\, {\mathcal C}_n\, $. For each $\, A_{\lambda }\, $ we will construct a homotopy of 
quasihomomorphisms, each representing an element of 
$\, E ( A_{\lambda }\, ,\, \Pi / \oplus\, A_{{\theta }_0} )\, $ (in principle one can use $KK$-theory but since  our stable projective limits are only well defined up to $E$-equivalence it seems more appropriate to pass to $E$-theory), connecting the evaluation map $\, A_{\lambda } \rightarrow \Pi / \oplus\, A_{{\theta }_0}\, $ with the trivial map and such that these homotopies constitute a coherent system of elements in $E$-theory with respect to inclusion. 
Let $\, \overline A_{\lambda }\, $ be the preimage of $\, A_{\lambda }\, $ in $\, \Pi\, C_n\, $ and choose a dense sequence $\, \{ x_k^{\lambda } \}\, $ of elements in the unit ball of $\, \overline A_{\lambda }\, $.
For each $\, n\, $ consider the image of  $\, \overline A_{\lambda }\, $ in 
$\, {\mathcal C}^{- \theta }_n \otimes A_{\theta }\, $ by the projection onto $\, {\mathcal C}_n\, $ composed with $\, {\delta }_{\theta }\, $ where $\, {\theta }_n \leq \theta \leq {\theta }_0\, $. The fibres of 
$\, {\mathcal C}_n^{- \theta }\, $ are indexed by the halfopen interval 
$\, ( {\theta }_n - \theta\, ,\, {\theta }_0 - \theta ]\, $ and there exists a small interval $\, I_n^{\lambda }\, $ around the origin such that replacing multiplication in the corresponding fibres by multiplication of the group $C^*$-algebra 
$\, C^* ( G )\, $ using the tautological shift $\, s_g^{\epsilon } \mapsto s_g^1\, $, the curvature (of the tautological shift) is uniformly small to the order of $\, 1 / n\, $ when restricted to the images of the finitely many elements $\, \{ x_1^{\lambda }\, ,\cdots ,\, x_n^{\lambda } \}\, $  for all $\,\theta\, $ in the compact interval $\, [ {\theta }_n\, ,\, {\theta }_0 ]\, $. Let 
$\, {\mathcal C}^{- \theta }_{n , +}\, $ be the image of $\, {\mathcal C}_n^{- \theta }\, $ corresponding to evaluation on the positive part $\, [ 0\, ,\, {\theta }_0 - \theta ]\, $ of the index set and similarly 
$\, {\mathcal C}_{n , -}^{- \theta }\, $ the image corresponding to evaluation on the negative part 
$\, ( {\theta }_n - \theta\, ,\, 0 ]\, $. The general idea is now to cut up a 
element $\, x\, $ of $\, {\mathcal C}_n^{- \theta }\, $ at the origin and to use the continuous analogue of
$\, j_{\sigma }\, $ to map $\, x_- \in {\mathcal C}_{n , -}^{- \theta}\, $ into 
$\, {\mathcal C}_{n , +}^{- \theta }\, $, then taking the difference of $\, x_+\, $ with the image of $\, x_-\, $ defines a quasihomomorphism into $\, {\mathcal C}_{n , +}^{-Ê\theta }\otimes A_{\theta }\, $.  The problem is that one wants the image of this construction to be trivial at the origin and the flip 
$\, j_0\, $ exchanges the two copies of $\, C^* ( G )\, $ in $\, A_0\, $ so one has to revert this process on some small interval around the origin. The idea is now to embed (and project)
$\, {\mathcal C}_n^{- \theta }\, $ to $\, \overline {\mathcal D}_n^{- \theta }\, $ where the latter denotes the semicontinuous field with fibres corresponding to the closed interval 
$\, [ \theta - {\theta }_0\, ,\, {\theta }_0 - \theta ]\, $ which is symmetric around the origin. Then embed 
$\, \overline {\mathcal D}_n^{- \theta }\, $ into 
$\, \overline {\mathcal D}_n^{- \theta } {\rtimes }_{\tau } {\mathbb Z}_2\, $ whose fibre at the origin is equal to $\, M_2 ( C^* ( G ) )\, $ so the flip becomes inner. Assuming stability one may choose a strictly continuous path of unitaries 
$\, \{ U_t \}\, $ in the stable multiplier algebra indexed by the interval 
$\, [ \theta  - {\theta }_0\, ,\, 0 ]\, $, each $\, U_t\, $ acting on the corresponding fibre at point $\, t\, $,
such that $\, U_0\, $ implements the flip $\,\tau\, $ and such that 
$\, U_t\, $ equals the identity at the negative end point of the interval $\, I_n^{\lambda }\, $ and all points beyond. Conjugating the image of $\, x\, $ in 
$\, \overline {\mathcal D}_n^{- \theta }\, $ with the unitary defined by the path $\, \{ U_t \}\, $ then taking the difference of the original image with the conjugated and flipped image defines a quasihomomorphism which evaluated at the origin is trivial. Tensoring with $\, A_{\theta }\, $ and evaluating at the point $\, {\theta }_0 - \theta\, $ (and its conjugate point) yields for each fixed 
$\,  {\theta }_n \leq \theta \leq {\theta }_0\, $ a quasihomomorphism from 
$\, {\mathcal C}_n\, $ to a certain subalgebra of 
$\, \bigl(\, ( A_{{\theta }_0 -\theta} \oplus A_{\theta - {\theta }_0}\,){\rtimes }_{\tau } {\mathbb Z}_2\,\bigr) \otimes A_{\theta }\, $. One notes that in case that $\, \theta\, $ is not contained in $\, I_n^{\lambda }\, $ 
the image lies in the subalgebra $\, A_{{\theta }_0 - \theta} \otimes A_{\theta }\, $. On the other hand, for 
$\, \theta \in I_n^{\lambda }\, $ the first factor is approximately equal $\, M_2 ( G )\, $ (or rather 
$\, M_2 ( C^* ( G ) ) \oplus M_2 ( C^* ( G ) )\, $ on symmetric fibres different from the origin) as far as the finitely many elements $\,\{ x_1^{\lambda }\, ,\cdots ,\, x_n^{\lambda } \}\, $ are concerned. Also the image will be contained in the subalgebra of this tensor product generated by diagonal elementary tensors of the form 
$\, \{ j^{\epsilon } s_g^{{\theta }_0 -\theta } \otimes j^{\epsilon } s_g^{\theta } \}\, $ (resp. also $\, \{ j^{\epsilon } s_g^{\theta - {\theta }_0} \otimes 
j^{\epsilon } s_g^{\theta } \}\, $ and $\,\tau\, $ if $\, \theta\in I_n^{\lambda }\, $) where $\, j\, $ is the element corresponding to the complex number $\, i\, $ in the underlying real Banach algebra of 
$\, {\mathcal A}^0\, $ and $\, \epsilon\, $ can take the values $\, 0\, $ or $\, 1\, $. Now for each pair of indices 
$\, ( {\theta }_0\, ,\, \theta )\, $ the $C^*$-algebra 
$\, \bigl(\, A_{{\theta }_0 - \theta }^0 \otimes A_{\theta }^0\, \bigr) {\otimes }_{\mathbb R} \mathbb C\, $ contains (commuting) copies of $\, A_{{\theta }_0 - \theta }\, $ and $\, A_{\theta }\, $ from which follows by the universal property of the maximal tensor product that one has a surjective homomorphism 
$$ A_{{\theta }_0 - \theta } \otimes A_{\theta } \twoheadrightarrow 
\bigl(\, A_{{\theta }_0 - \theta }^0 \otimes A_{\theta }^0\,\bigr) {\otimes }_{\mathbb R} \mathbb C\> . $$
The subalgebra of the right side generated by diagonal elementary tensors 
$\, \{ s_g^{{\theta }_0 - \theta } \otimes s_g^{\theta } \}\, $ is isomorphic with 
$\, A_{{\thetaÊ}_0}\, $.  To see this note that the tensor product algebras 
$\, A_{\sigma }^0 \otimes A_{\theta }^0\, $ arise as fibres of the semicontinuous field $\, C^* ( \overline G \times \overline G ) \simeq C^* ( \overline G ) \otimes C^* ( \overline G )\, $ over $\,\widehat K \times \widehat K\, $. Restricting the semicontinuous field to the $C^*$-subalgebra $\, C^* ( \delta ( \overline G ) \times K )\, $ where $\,\delta ( \overline G )\, $ denotes the diagonal in $\, \overline G \times \overline G\, $ and the other factor $\, K\, $ is identified with the central subgroup $\, K\, $ in the second copy of $\,\overline G\, $, one finds that (the injective image in $\, C^* ( \overline G \times \overline G )\, $ of)
$\, C^* ( \delta ( \overline G ) \times K )\, $ is isomorphic with the tensor product 
$\, C^* ( \overline G ) \otimes C^* ( K )\, $ so that each of its fibres corresponds to some fibre 
$\, A_{{\theta }_0 }^0 \, $ of $\, C^* ( \overline G )\, $, and on the other hand it is sent to the diagonal in 
$\, A_{{\theta }_0 - \theta }^0 \otimes A_{\theta }^0\, $ modulo the ideal of functions vanishing at the point 
$\, ( {\theta }_0  - \theta , \theta )\, $, so that both 
$C^*$-algebras must be isomorphic. For supposing that the map on fibres is not injective consider an element in its kernel and, choosing some small interval around the given point where the tautological shift is approximately contractive for the given element (resp. a suitable approximation from the dense subalgebra of finite sums in the generators $\,\{ s_g^{{\thetaÊ}_0}Ê\}\, $), also the images in the tensor product of fibres corresponding to points in the interval must become uniformly small for the shifted elements by upper semicontinuity if the interval is small enough. Then multiplying by some $\delta $-shaped positive realvalued function of norm one whose domain is contained inside the interval and equal to 1 at 
$\, {\theta }_0\, $ one gets a contradiction because the norm of the image of the constructed extension is arbitrarily small, while on the other hand the injection $\, C^* ( \delta ( \overline G ) \times K ) \hookrightarrow C^* ( \overline G \times \overline G )\, $ is an isometry. 
Returning to the general argument one  
combines the quasihomomorphisms constructed above for all values of 
$\, \theta\in [ {\theta }_n\, ,\, {\theta }_0 ]\, $ and for all $\, n\in \mathbb N\, $. This yields a quasihomomorphism 
from $\, \overline A_{\lambda }\, $ to the outer direct product 
$\, \Pi / \oplus\, \bigl( M_2 ( A_{{\theta }_0} ) \otimes C [ {\theta }_n\, ,\, {\theta }_0 ] \bigr)\, $. Restriction of this quasihomomorphism to the direct sum 
$\, \oplus\, {\mathcal C}_n\, $ is trivial, so that by standard techniques the construction yields an element of 
$$ E \bigl( A_{\lambda }\, ,\, \Pi / \oplus\, ( A_{{\theta }_0} \otimes C [ {\theta }_n\, ,\, {\theta }_0 ] ) \bigr) $$
which when evaluated at the end points $\, {\theta }_0\, $ gives the trivial element while when evaluated at the sequence of points 
$\, \{ {\theta }_n \}\, $ it gives the element induced by the evaluation map 
$\, A_{\lambda } \rightarrow \Pi / \oplus\, A_{{\theta }_0}\, $. But 
$\, \Pi / \oplus\, ( A_{{\theta }_0} \otimes C [ {\theta }_n\, ,\, {\theta }_0 ] )\, $ is $E$-equivalent to 
$\, \Pi / \oplus\, A_{{\theta }_0}\, $ by either projection, which follows from the corresponding result for the direct sum and the direct product using the six-term exact sequence of $E$-theory. Then 
the evaluation map $\, \Pi / \oplus\, {\mathcal C}_n \twoheadrightarrow \Pi / \oplus\, A_{{\theta }_0}\, $ is trivial in $K$-theory as asserted. The rest of the argument is a simple diagram chase.
Suppose that for some $\, x \in K_1 ( C_0 ( \mathbb R ) )\, $ the composition 
$$ K_1 ( C_0 ( \mathbb R ) ) \longrightarrow K_1 ( \mathcal S_0 )\longrightarrow K_1 ( \mathcal S_n ) $$ 
is trivial for every $\, n > 0\, $. Then the image $\, y\, $ in $\, K_1 ( \mathcal S_\infty )\, $ is contained in the subgroup 
$\, \lim^1 K_0 ( \mathcal S_n )\, $ which can be identified with the cokernel of 
$$ \Pi / \oplus\, K_1 ( {\mathcal S}_n ) \buildrel \Pi / \oplus\, ( id - {\rho }_n )_*\over\largerightarrow 
\Pi / \oplus\, K_1 ( {\mathcal S}_n ) \> . $$
Choose a representative $\, \xi\, $ for this element in $\, K_1 ( \Pi / \oplus\, {\mathcal S}_n )\, $. Its image in $\, K_1 ( \Pi / \oplus\, {\mathcal C}_n )\, $ has a preimage for $\, \Pi / \oplus\, ( id - {\rho }_n )_*\, $ which maps to zero in $\, K_1 ( \Pi / \oplus\, A_{{\theta }_0} )\, $ by evaluation at $\, {\theta }_0\, $, hence it lifts to an element of $\, K_1 ( \Pi / \oplus\, {\mathcal O}_n )\, $ which again maps to an element $\, \zeta\, $ of
$\, K_1 ( \Pi / \oplus\, {\mathcal S}_n )\, $. The difference $\, \xi - \Pi / \oplus\, ( id - {\rho }_n )_* ( \zeta )\, $ represents the same element of the cokernel as $\, \xi\, $, and lifts to 
$\, K_1 ( \Pi / \oplus\, {\mathcal S}_0 )\, $. But the latter group has trivial cokernel for 
$\, \Pi / \oplus\, ( id - {\rho }_n )_*\, $ so that also $\, y\, $ must be trivial. Then the image of $\, x\, $ in 
$\, K_1 ( {\mathcal S}_0 )\, $ lifts to an element of $\, K_0 ( {\mathcal C}_{\infty } )\, $ which maps to the zero element of $\, K_0 ( A_{{\theta }_0} )\, $ under evaluation, so it lifts again to an element of 
$\, K_0 ( {\mathcal O}_{\infty } )\, $ whose image in $\, K_1 ( {\mathcal S}_0 )\, $ is trivial by exactness, giving a contradiction. 
\par\smallskip\noindent
In the general case (of nonfinitely generated $\, K\, $) one proceeds similarly. 
Consider the push out of the central extension $\, \overline G\, $ by the map 
$\, K \rightarrow K \otimes \mathbb Q\, $ and check that it is well defined from the fact that $\, K\, $ is central. If the map of the proposition is not rationally injective there exists a nontorsion element of 
$\, K \, $ in its kernel which may be completed to a vector space basis of $\, K \otimes \mathbb Q\, $. As above one now divides by the complementary subspace  as to obtain a central extension of 
$\, G\, $ by $\, \mathbb Q\, $ which we continue to denote $\,\overline G\, $, and its kernel by $\, K\, $.
The kernel $\, \mathcal J\, $ of the surjection $\, C^* ( \overline G ) \twoheadrightarrow C^* ( G )\, $ is again a module over the central subalgebra $\, C^*_0 ( K )\, $ isomorphic to the algebra of continuous functions on the compact dual $\,\widehat K\, $ of $\, K = \mathbb Q\, $ vanishing at the identity element 
$\, \{ 1 \}\, $, 
which is the maximal $C^*$-completion of the algebra of finite sums as above with coefficients in 
$\, C^*_0 ( K )\, $. The algebra $\, C^* ( K )\, $ can be realized as continuous functions on 
the real line which are periodic with integer period $\, q\, $ for some $\, q \geq 1\, $. 
$\, C^*_0 ( K )\, $ contains those functions vanishing at the point $\, 0\, $. Instead of a single multiplier 
$\, u\, $ defining multiplication on the algebra of finite sums 
$\, \{ \sum_g\, a_g s_g \}\, $ with respect to some given section $\, G \nearrow \overline G\, $ one now has to consider the family of functions $\,\{ u_q\,\vert\, q \in \mathbb N \}\, $ defined by 
$\, u_q ( t ) = e^{2 \pi i t / q}\, $, i.e. for each $\, g , h \in G\, $ there exist natural numbers 
$\, p = p ( g , h )\, ,\, q = q ( g , h )\, $ such that the relation 
$\, s_g\, s_h = u_q^p\, s_{gh}\, $ holds in $\, C^* ( \overline G )\, $. 
Let $\, C^*_{0 , q} ( K )\, $ denote the ideal in $\, C^*_0 ( K )\, $ of functions vanishing at all integral multiples of $\, q\, $. Then $\, C^*_0 ( K )\, $ is the inductive limit of the net 
$\, \{ C^*_{0 , q} ( K ) \}\, $ directed by the natural inclusions. Let $\, \{ {\mathcal J}_q \}\, $ be the corresponding increasing net of ideals with limit $\,\mathcal J\, $. In order to show injectivity of the restriction $\, K_1 ( C^*_{0 , q} ( K ) ) \rightarrowtail K_1 ( {\mathcal J}_q )\, $ one proceeds exactly as above, deforming the action of 
$\, u_q\, $ by some parameter $\, \theta \in [ 0 , 1 ]\, $ such that $\, u_{q , \theta }\, $ is periodic of period 
$\, q\, $ and equals the identity function on the interval $\, [ 0 , q \theta ]\, $ whereas on the interval 
$\, [ q \theta , q ]\, $ one has $\, u_{q , \theta } ( t ) = u ( t - q \theta )\, $. Also this deformation can be extended in a compatible way to any multiplier $\, u_{mq}\, $ on taking the m-th root of the deformation of 
$\, u_q\, $. Although $\, u_{mq , \theta }\, $ is not continuous it does define a multiplier of 
$\, C^*_{0 , q} ( K )\, $, and $\, u_{q ,1} \equiv 1\, $. The  argument is now the same as for finitely generated $\, K\, $, considering the corresponding continuous field of $C^*$-algebras $\, \{ {\mathcal J}_{q , \theta } {\} }_{\theta \in [ 0 , 1 ]}\, $. Passing to inductive limits gives the general result.
This completes the proof\qed  
\par\bigskip\noindent
{\it Remark.}\quad If $\, G = \mathbb Z \oplus \mathbb Z\, $ the fibres $\, A_{\theta }\, $ of the group algebra of the universal central extension as above are just the irrational rotation algebras which have been studied extensively by various authors (cf. \cite{Rf2} , \cite {Rf3} , \cite{C}). It is well known that in this particular case the fibres are all $KK$-equivalent so the cone-like extensions corresponding to the $\, {\mathcal C}_n\, $ are (globally) contractible which is a much stronger statement than the local property derived above. For a general group $\, G\, $ there is an interesting characterization of the stable fibres 
$\, \mathcal K \otimes A_{\theta }\, $ as a crossed product of the algebra of compact operators on 
$\, \mathit l^2 ( G )\, $ by an outer group action of $\, G\, $ which goes as follows. 
Consider the left regular representation of $\, G\, $ on 
$\, \mathit l^2 ( G )\, $. If $\, \{ {\epsilon }_k\, \vert\, k \in G \}\, $ denotes the standard orthonormal basis, then the left regular representation of $\, A_{\theta }\, $ (compare \cite{Rf}) is defined by 
$$ \lambda ( s_g^{\theta } ) ( {\epsilon }_k ) = exp ( 2\pi i ( 1 - {\theta } ) n ( g , k ) )\cdot {\epsilon }_{gk}  \> . $$
The enveloping $C^*$-algebra $\,\lambda ( A_{\theta } )\, $ acts on the commutative $C^*$-algebra 
$\, C_0 ( G )\, $ embedded in the usual way by diagonal matrices in $\,\mathcal K ( \mathit l^2 ( G ))\, $ by 
$\, \lambda ( s_g^{\theta })\, p_h\, \lambda ( s_g^{\theta } )^* = p_{gh}\, $ where $\, p_h\, $ denotes the onedimensional orthogonal projection onto the subspace spanned by $\, {\epsilon }_h\, $. The corresponding "crossed product" $\, C_0 ( G ) \rtimes \lambda ( A_{\theta } )\, $ is contained in  
$\,\mathcal K ( \mathit l^2 ( G ))\, $. Since the elements $\, \{ \lambda ( s_g^{\theta } ) \}\, $ differ from 
the generators $\, \{ \lambda ( s_g^1 ) \}\, $ of $\, \lambda ( G )\, $ in the left regular representation only modulo diagonal invertible matrices which are multipliers of $\, C_0 ( G )\, $, the twisted crossed product is in fact equal to $\,\mathcal K ( \mathit l^2 ( G ) )\, $.
The integer-valued function $\, n ( g , h )\, $ is defined by the chosen section 
$\, G \nearrow\overline G\, $. Changing from one section to another amounts to multiplying each generator 
$\, s_g\, ,\, g \neq 1\, $ with 
a complex number $\, e^{2 \pi i \theta m ( g )}\, $ where $\, m : G\backslash \{ 1 \} \rightarrow \mathbb Z\, $ is some integer valued (or more generally real valued) function. This changes the 2-cocycle 
$\, n\, $ to $\, \widetilde n ( g , h ) = n ( g , h ) + m ( g ) + m ( h ) - m ( gh )\, $. It is clear that one can always find a selfadjoint section satisfying $\, n ( g , g^{-1} ) = 0\, $ for all $\, g \in G\, $ (if $\, G\, $ has $2$-torsion one has to consider halfinteger-valued modifications $\, m\, $). This implies that 
$\, s_{g^{-1}} = s_g^*\, $ and hence $\, n ( g , h ) = - n ( h^{-1} , g^{-1} )\, $. For $\, 0 < \theta \leq 1\, $ put 
$\, {\theta }' = 1 - \theta\, ,\, \omega = exp ( 2\pi i {\theta }')\, ,\, {\omega }' = exp ( 2\pi i \theta )\, $.
Put $\, u_g^{\theta } = \lambda ( s_g^{{\theta }'} )\, $ and let $\, G\, $ act on the $C^*$-algebra 
$\, C_0 ( G ) \rtimes \lambda ( A_{{\theta }'} )\, $ by the formula 
$\, {\alpha }^{\theta }_g (  p_h ) = p_{gh}\, $ plus the following twisted adjoint operation on the image of 
$\, A_{{\theta }'}\, $
$$ {\alpha }^{\theta }_g ( u_h^{\theta } ) \> =\> 
{\omega }^{n ( g , h^{-1} ) - n ( g , hg^{-1} )}\cdot u_{ghg^{-1}}^{\theta } $$
and check that this defines an action of $\, G\, $ compatible with the actions of $\, G\, $ and 
$\, A_{{\theta }'}\, $ on $\, C_0 ( G )\, $ so that the action extends to an action 
$\, {\alpha }^{\theta }\, $ on $\, \mathcal K ( \mathit l^2 ( G ))\, $. The corresponding (maximal) crossed product is not (in an obvious way) isomorphic to $\, \mathcal K ( \mathit l^2 ( G )) \otimes C^* ( G )\, $ since, although any automorphism of $\,\mathcal K\, $ is inner, the unitaries implementing the action do not define a (unitary) representation of $\, G\, $. Let 
$\, \{ v^{\theta }_g\,\vert\, g \in G \}\, $ denote the unitaries corresponding to group elements in the multipliers of the  crossed product algebra. Putting 
$\, w^{\theta }_g = u^{\theta }_{g^{-1}}\cdot v^{\theta }_g\, $ defines a representation of $\, A_{\theta }\, $ mapping $\, s^{\theta }_g\, $ to the unitary $\, w^{\theta }_g\, $. Also the  
$\, \{ w^{\theta }_g \}\, $ commute with $\, \mathcal K ( \mathit l^2 ( G ))\, $.
Thus one gets an (isomorphic) map $\, \mathcal K ( \mathit l^2 ( G )) \otimes A_{\theta } \longrightarrow 
\mathcal K ( \mathit l^2 ( G )) \rtimes_{{\alpha }^{\theta }} G\, $ for each $\,\theta\, $. This result should remain valid on replacing the maximal group algebras and maximal crossed products by the corresponding reduced ones throughout.
\par

\end{document}